%%%%%%%%%%%%%%%%%%%%%%%%%%%%%%%%%%%%%%%%%%%%%%%%%%%%%%%%%%%%%%%%%
%  Version 1.0 started on 10/14/00, finished on 10/28/00
%%%%%%%%%%%%%%%%%%%%%%%%%%%%%%%%%%%%%%%%%%%%%%%%%%%%%%%%%%%%%%%%%

\input amstex
\documentstyle{amsppt}

\hcorrection{19mm}

\nologo
\NoBlackBoxes

%\NoRunningHeads
%\mag=1200
%\hsize=31 pc
%\vsize=44 pc
%\baselineskip 15pt
%\hcorrection{5mm}

\topmatter
\title  Completely tubing compressible tangles and standard graphs in
genus one 3-manifolds
\endtitle
\author  Ying-Qing Wu$^1$
\endauthor
\leftheadtext{Ying-Qing Wu}
\rightheadtext{Completely tubing compressible tangles}
\address Department of Mathematics, University of Iowa, Iowa City, IA
52242
\endaddress
\email  wu\@math.uiowa.edu
\endemail

\thanks  $^1$ Partially supported by NSF grant \#DMS 9802558
\endthanks

\abstract We prove a conjecture of Menasco and Zhang that if a tangle
is completely tubing compressible then it consists of at most two
families of parallel strands.  This is related to problems of graphs
in 3-manifold.  A 1-vertex graph $\Gamma$ in a 3-manifold $M$ with a
genus 1 Heegaard splitting is standard if it consists of one or two
parallel sets of core curves lying in the Heegaard splitting solid
tori of $M$ in the standard way.  The above conjecture then follows
from the theorem which says that a 1-vertex graph in $M$ is standard
if and only if the exteriors of all its nontrivial subgraphs are
handlebodies.  \endabstract

\endtopmatter

\document

\redefine\Gamma{G}
\define\proof{\demo{Proof}}
\define\endproof{\qed \enddemo}
\define\a{\alpha}
\redefine\b{\beta}

\define\r{\gamma}

\redefine\bdd{\partial}
\define\Int{\text{\rm Int}}

In this paper, a tangle is a pair $(W, t)$, where $W$ is a compact
orientable 3-manifold with $\bdd W$ a sphere, and $t = \a_1 \cup
... \cup \a_n$ a set of mutually disjoint properly embedded arcs in
$W$, called the strands.  We denote by $N(t)$ a regular neighborhood
of $t$, and by $\eta(t)$ an open neighborhood of $t$, i.e.\ $\eta(t) =
\Int N(t)$.  Denote by $X = X(t)$ the tangle space $W - \eta(t)$, and
by $P$ the planar surface $\bdd W \cap X = \bdd W - \eta(\bdd t)$.
Let $A_i$ be the annulus $\bdd N(\a_i) \cap X$.  Thus $\bdd X = P \cup
(\cup A_i)$.

Following Gordon [G], we say that a set of curves $\{c_1, ..., c_k\}$
on the boundary of a handlebody $H$ is {\it primitive\/} if there
exist disjoint disks $D_1, ..., D_k$ in $H$ such that $\bdd D_i$
intersects $\cup c_j$ transversely at a single point lying on $c_i$.
A set of annuli is primitive if their core curves form a primitive
set.

Denote by $F_i = P \cup A_i$, call it the $A_i$-tubing surface of $P$.
The surface $P$ is $A_i$-tubing compressible if $F_i$ is compressible,
and it is completely $A_i$-tubing compressible if $F_i$ can be
compressed until it becomes a set of annuli parallel to $\cup _{j\neq
i} A_j$.  Equivalently, $P$ is completely $A_i$-tubing compressible if
$X$ is a handlebody, and the set of annuli $\cup _{j\neq i} A_j$ is
primitive on $\bdd X$.  The tangle $(W,t)$ is {\it completely tubing
compressible\/} if it is completely $A_i$-tubing compressible for all
$i$.  Such tangles arise naturally in the study of reducible surgery
on knots.  It has been shown in [CGLS] that if some surgery on a
hyperbolic knot $K$ produces a nonprime manifold $M$, then either the
knot complement contains a closed essential surface, or there is a
reducing sphere $S$ cutting $(M, K')$ into two non-split completely
tubing compressible tangles, where $K'$ is the core of the Dehn
filling solid torus.

Define a {\it band\/} in $W$ to be an embedded disk $D$ in $W$ such
that $D \cap \bdd W$ consists of two arcs on $\bdd D$.  A
subcollection of strands $t' = \{\a_1, ... ,\a_k\}$ of $t$ is {\it
parallel\/} if there is a band $D$ such that $D\cap t = t'$.

Define a {\it core arc\/} to be an arc $\a$ in $W$ such that $W -
\eta(\a)$ is a solid torus.  Because of uniqueness of Heegaard
splittings of $S^3$, $S^2 \times S^1$ and lens spaces, it is easy to
see that $W$ has at most two core arcs up to isotopy, and one if $W$
is a punctured $S^3$, $S^2\times I$ or $L(p,1)$.  However, a set of
core arcs may contain arbitrarily many parallel families.  This is the
same phenomenon as links in $S^3$: A link $L$ of $n$ components may
have the property that all of its components are trivial knot (so the
components are isotopic to each other in $S^3$), but the components of
$L$ are mutually non-parallel in the sense that they do not bound an
annulus with interior disjoint from the link.  The following theorem
proves a conjecture of Menasco and Zhang [MZ, Conjecture 5], which
shows that this phenomenon will not happen if $(W, t)$ is a completely
tubing compressible tangle.  I would like to thank Menasco and Zhang
for posting the conjecture.

\proclaim{Theorem 1} If $(W, t)$ is a completely tubing compressible
tangle, then $t$ consists of at most two families of parallel core
arcs.  \endproclaim

The problem is related to graphs in 3-manifolds.  Let $M = \hat W$ be
the union of $W$ and a 3-ball $B$, and let $\Gamma = \hat t$ be the
union of $t$ and the straight arcs in $B$ connecting $\bdd t$ to the
central point $v$ of $B$.  Thus we have a graph $\hat t$ in the closed
3-manifold $\hat W$ with one vertex $v$ and $n$ edges $e_1, ..., e_n$
corresponding to the arcs $\a_1, ..., \a_n$ of $t$.  A graph $\Gamma$
is {\it nontrivial\/} if it contains at least one edge.  The {\it
exterior\/} of a graph $\Gamma$ in a 3-manifold $M$ is $E(\Gamma) = M
- \eta(\Gamma)$.  The following lemma translated the completely tubing
compressible condition to a condition about $\hat t$ in $\hat M$.

\proclaim{Lemma 2} The tangle $(W, t)$ is completely tubing
compressible if and only if the exterior of any nontrivial subgraph of 
$\hat t$ in $\hat W$ is a handlebody.  \endproclaim

\proof Let $A_i$ be the annulus $\bdd N(\a_i) \cap \bdd X$.  The
exterior of a subgraph $\Gamma'$ of $\Gamma = \hat t$ in $\hat W$ is
the same as the exterior of the corresponding strands of $t$ in $W$,
which can be obtained from $X = W - \eta(t)$ by attaching 2-handles to
those annuli $A_i$ corresponding to the edges $e_i$ in $\Gamma -
\Gamma'$.  Therefore the condition that the exterior of any nontrivial
subgraph of $\hat t$ in $\hat W$ is a handlebody implies that
attaching 2-handles to $X$ along any proper subset of $\cup A_i$
yields a handlebody.  By [G, Theorem 1] this implies that any proper
subset of $\cup A_i$ is a primitive set on $\bdd X$.  Hence $(W, t)$
is completely tubing compressible.

On the other hand, if $(W, t)$ is completely tubing compressible, and
$\Gamma'$ is a proper subgraph of $\Gamma$ which does not contain the
edge $e_i$, say, then the set $\cup _{j\neq i} A_j$ is primitive on
$\bdd X$, and since the exterior $E(\Gamma')$ of $\Gamma'$ can be
obtained by attaching 2-handles to $X$ along a subset of primitive set
$\cup _{j\neq i} A_j$, it follows that $E(\Gamma')$ is a handlebody.
\endproof

The classification problem for completely tubing compressible tangles
now becomes a classification problem for 1-vertex graphs $\Gamma$ in a
3-manifold $M$ which have the property that the exteriors of all its
nontrivial subgraphs are handlebodies.  Since the exterior of a
regular neighborhood of an edge of $\Gamma$ is a solid torus, $M$ has
a Heegaard splitting of genus 1, hence it must be $S^3$, $S^2\times
S^1$, or a lens space $L(p,q)$.  Since $L(p,q) \cong L(p, -q) \cong
L(p, p-q)$ up to (possibly orientation reversing) homeomorphism, we
may always assume that $1\leq q \leq p/2$.  When $M$ is $S^3$, it
follows from [G, Theorem 1] that the complement of any subgraph of
$\Gamma$ is a handlebody if and only if $\Gamma$ is planar, i.e., it
is contained in a disk in $S^3$.  Scharlemann and Thompson [ST]
generalizes this to all abstractly planar graphs in $S^3$.  See also
[Wu2] for an alternative proof.  For the general case, we need the
following definitions.

A {\it $v$-disk\/} $D$ in $M$ is the image of a map $f: D^2 \to M$
such that $f$ is an embedding except that it identifies two boundary
points of $D^2$ to a point $v$ in $M$.  The boundary of $D$ is $\bdd D
= f(\bdd D^2)$.  A $v$-disk $D$ in a solid torus $V$ is {\it
standard\/} if (i) $D \cap \bdd V = v$, and (ii) $D$ is rel $v$
isotopic to a $v$-disk $D'$ on $\bdd V$, which is longitudinal in the
sense that there is a meridional disk $\Delta$ of $V$ such that $D\cap
\Delta$ is a nonseparating arc on $D$.  We remark that it is important
to require that the above isotopy be relative to $v$ as that
guarantees that the exterior of $D$ is a handlebody.

A graph with a single vertex is called a {\it 1-vertex graph}.  Such a
graph is connected, and all of its edges are loops.  A 1-vertex graph
$\Gamma = e_1 \cup ... \cup e_k$ in $V$ with vertex $v$ is {\it in
standard position\/} if it is contained in a standard $v$-disk $D$ in
$V$.  In this case we also say that the edges of $\Gamma$ are
parallel.

Let $V_1 \cup V_2$ be a genus one Heegaard splitting of a closed
3-manifold $M$.  Then a 1-vertex graph $\Gamma$ in $M$ is {\it in
standard position\/} (relative to the Heegaard splitting) if either
(i) $M$ is homeomorphic to $S^3$, $S^2\times S^1$ or $L(p,1)$, and
$\Gamma$ is contained in a single standard $v$-disk in $V_1$ or $V_2$,
or (ii) $M$ is homeomorphic to $L(p,q)$ with $2\leq q < p/2$, and
$\Gamma$ is contained in two standard $v$-disks, one in each $V_i$.  A
1-vertex graph $\Gamma$ in $M$ is {\it standard\/} if it is isotopic
to a graph in standard position.  Since genus one Heegaard splittings
of 3-manifolds are unique up to isotopy [W, BO, S], this is
independent of the choice of $(V_1,V_2)$.  The following theorem
characterizes standard graphs in 3-manifolds.

\proclaim{Theorem 3} A nontrivial 1-vertex graph $\Gamma$ in a closed
orientable 3-manifold $M$ is standard if and only if the exterior of
any nontrivial subgraph of $\Gamma$ is a handlebody.  \endproclaim

It should be noticed that the 3-manifold $M$ in the theorem must be
$S^3$, $S^2\times S^1$, or a lens space.  For if $\Gamma$ is standard
then by definition $M$ has a genus one Heegaard splitting.  On the
other hand, if the exterior of any nontrivial subgraph of $\Gamma$ is
a handlebody, then in particular the exterior of an edge of $\Gamma$
is a solid torus, so again $M$ has a genus one Heegaard splitting.
Therefore $M$ must be one of the above manifolds.

The following lemma proves the easy direction of the theorem.

\proclaim{Lemma 4} If a 1-vertex graph $\Gamma$ in a 3-manifold $M$ is
standard, then the exterior of any nontrivial subgraph $\Gamma'$ of
$\Gamma$ is a handlebody.  \endproclaim

\proof Clearly a subgraph of $\Gamma$ is still standard, hence we need
only prove the lemma for $\Gamma' = \Gamma$.  Let $(V_1, V_2)$ be a
genus one Heegaard splitting of $M$, and assume that $\Gamma$ is
contained in the union of $D_1 \cup D_2$, where $D_i$ is a standard
$v$-disk in $V_i$.  (The case that $\Gamma$ is contained in a single
standard $v$-disk is similar and simpler.)  Put $\Gamma_1 = \Gamma
\cap D_1 = e_1 \cup ... \cup e_{r-1}$ and $\Gamma_2 = \Gamma \cap
D_2 = e_r \cup ... \cup e_n$.

From definition one can see that the manifold $V_i - \eta(D_i)$ is a
product $F_i \times I$, where $F_i$ is a once punctured torus.
Therefore $X = M - \eta (D_1\cup D_2)$ is still a product of $I$ and
a once punctured torus, which is a handlebody.  One can choose a
regular neighborhood $N(\Gamma)$ of $\Gamma$ in $M$ so that it is
contained in $N(D_1\cup D_2)$, and the closure of each component of
$N(D_1\cup D_2) - N(\Gamma)$ is a 3-ball $H_i$ intersecting $\bdd
N(D_1\cup D_2)$ at two disks.  Now $M - \eta(\Gamma)$ is the union of
$X$ and the $H_i$.  Since each $H_i$ can be considered as a 1-handle
attached to $X$, it follows that $M - \eta(\Gamma)$ is a handlebody.
\endproof

The following lemma proves the other direction of
Theorem 3 under an extra assumption, which by [MZ, Lemma 1] implies
that $M = S^3$ or $S^2 \times S^1$.  

\proclaim{Lemma 5} Let $\Gamma$ be a 1-vertex graph in a closed
orientable 3-manifold such that the exterior of any nontrivial
subgraph of $\Gamma$ is a handlebody.  Let $W = M - \eta(v)$, and $X =
M - \eta(\Gamma)$.  If $P = \bdd W \cap X$ is compressible, then
$\Gamma$ is standard.  \endproclaim

\proof Let $D$ be a compressing disk of $P$.  First assume that $D$ is
separating in $W$, cutting $W$ into $W_1$ and $W_2$.  Let $\Gamma_i$
be the subgraph of $\Gamma$ consisting of edges whose intersection
with $W$ is contained in $W_i$.  Each $\Gamma_i$ is nontrivial as
otherwise $\bdd D$ would be trivial on $P$, contradicting the fact
that it is a compressing disk.  Now $W_i$ is contained in the exterior
of $\Gamma_j$ ($j\neq i$), which by assumption is a handlebody.  Since
$\bdd W_i = S^2$ and handlebodies are irreducible, it follows that
$W_i$ are 3-balls, hence $W$ is also a 3-ball, so $M = S^3$.  In this
case by [G, Theorem 1] or [ST], the graph $\Gamma$ is planar in
$S^3$, which is easily seen to be equivalent to the condition that it
is standard.

Now assume the $D$ is non-separating in $W$.  In this case $W$ cannot
be a 3-ball or punctured lens space, so it must be a punctured
$S^2\times S^1$, and $D$ cuts $W$ into $W' = S^2 \times I$.  The
manifold $X' = W' - \eta(t)$ is obtained from $X$ by cutting along a
nonseparating disk $D$, so it is a handlebody of genus $n-1$, and
attaching 2-handles to any proper subset of $\cup A_i$ yields a
handlebody.  By [G, Theorem 2], the set $\cup A_i$ is standard on
$\bdd X'$, which implies that there is a band $D' = C \times I$ in $W'
= S^2 \times I$ containing $t = \Gamma \cap W$.  It is clear that such
a band $D$ extends to a standard $v$-disk $D''$ in $M = S^2\times S^1$
containing $\Gamma$.  \endproof

A trivial arc in a solid torus $V$ is one which is rel $\bdd$ isotopic
to an arc on $\bdd V$.  Given a $(p,q)$ curve $\gamma$ on $\bdd V$
(running $p$ times along the longitude) and a trivial arc $\a$ in $V$
disjoint from $\gamma$, the {\it jumping number\/} of $\a$ relative to
$\gamma$, denoted by $j(\a, \r)$, is defined as the minimal
intersection number between $\gamma$ and all arcs on $\bdd V$ which is
rel $\bdd$ isotopic to $\a$.  Clearly we have $0\leq j(\a,\r) \leq
p/2$, and the arc on $\bdd V$ which is isotopic to $\a$ and
intersects $\r$ at $j(\a,\r)$ points must intersect $\r$ always in
the same direction.  The following lemma is essentially [MZ,
Proposition 3].  The proof here is more straight forward.

\proclaim{Lemma 6} Let $\r$ be a $(p,q)$ curve on the boundary $T$ of
a solid torus $V_1$ with $1\leq q \leq p/2$.  Let $\a$ be a trivial arc in
$V_1$ with boundary disjoint from $\r$, and let $\b$ be an arc on $T$
disjoint from $\r$, connecting the two endpoints of $\a$.  Let
$L(p,q)$ be the lens space obtained by gluing a solid torus $V_2$ to
$V_1$ such that $\r$ bounds a meridional disk in $V_2$.  If the exterior of
$\a\cup \b$ in $L(p,q)$ is a solid torus, then the jumping number
$j(\a, \r)$ of $\a$ relative to $\r$ is either $1$ or $q$.
\endproclaim

\proof Since $\pi_1 L(p,q) = {\Bbb Z}_p$, by choosing an orientation
properly every curve $\delta$ in $L(p,q)$ represents a unique element
$[\delta]$ between $0$ and $p/2$.  By definition $\a$ is isotopic rel
$\bdd$ to an arc $\a'$ on $T$ intersecting $\r$ transversely at
$j(\a,\r)$ points in the same direction.  Thus if we choose the core
curve of $V_2$ as a generator of $\pi_1 L(p,q) = \Bbb Z_p$, then the
curve $\delta = \a \cup \b$ represents the number $j(\a, \r)$ in $\Bbb
Z_p$.  On the other hand, since the exterior of $\delta$ is a solid
torus, by uniqueness of Heegaard splittings of lens spaces [BO], the
curve $\delta$ is isotopic to the core of either $V_1$ or $V_2$, which
represents the elements $1$ and $q$ in $\Bbb Z_p$, respectively, hence
the result follows.  \endproof

\proclaim{Lemma 7} Theorem 3 is true if $M = L(p,q)$ and $\Gamma$ has
at most two edges. \endproclaim

\proof If $\Gamma$ has only one edge $e_1$, then $V_1=N(e_1)$ and $V_2
= M - \Int V_1$ form a genus one Heegaard splitting of $L(p,q)$.  By
an isotopy we may deform $e_1$ to standard position in $V_1$, and the
result follows.

We now assume that $\Gamma = e_1 \cup e_2$.  Let $V_1 = N(e_1)$, and
$V_2 = M - \Int V_1$, which by assumption is a solid torus.  Since
$e_2$ intersects $e_1$ at the vertex $v$ of $\Gamma$, we may assume
that $e_2\cap V_1$ is an unknotted arc lying on a meridional disk $D'$
of $V_1$.  Let $D$ be another meridional disk of $V_1$ disjoint from
$D'$, and let $\r$ be the curve $\bdd D$ on $T = \bdd V_i$.  Since $M$
is a lens space $L(p,q)$, $\r$ is a $(p,q)$ curve on $T$ with respect
to some longitude-meridian pair of $V_2$.  Let $\a$ be the embedded
arc $e_2 \cap V_2$ in $V_2$.  The boundary of $\a$ lies on the curve
$\r' = \bdd D'$, which is a parallel copy of $\r$.

Note that $V_2 - \eta(\a) = M - \eta(\Gamma)$, so by assumption it is
a handlebody, denoted by $H$.  The frontier of $N(\a)$ is an annulus
$A$ which must be primitive on $H$ because when attaching the 2-handle
$N(\a)$ to $H$ along $A$ we obtain the solid torus $V_2$.  It follows
that the core curve $\a$ of the attached 2-handle $N(\a)$ is a trivial
arc in $V_2$.

Let $\b$ be an arc on $\bdd \r'$ connecting the two endpoints of $\a$.
Then $\b$ is isotopic to the arc $e_2\cap V_1$ on the disk $D'$, hence
the curve $\a \cup \b$ is isotopic to $e_2$, which by assumption has
exterior a solid torus in $L(p,q)$.  Therefore by Lemma 6 the jumping
number $j(\a, \r)$ is either $1$ or $q$.  By definition $\a$ is
isotopic rel $\bdd$ to an arc $\a'$ on $T$ intersecting $\r$
transversely at $j(\a,\r)$ points in the same direction.

First assume that $j(\a,\r) = 1$.  Then $e_2' = \a' \cup \b$ is a
simple closed curve on $T$ intersecting the meridian curve $\r$ of
$V_1$ transversely at a single point, hence it is a longitude of
$V_1$.  Since $\b$ lies on $\bdd D'$ and $e_2\cap V_1$ is an arc on
$D'$, there is an isotopy of $\Gamma \cap V_1$ in $V_1$ such that $e_2
\cap V_1$ is deformed to the arc $\b$, and $e_1$ to a loop $e'_1$ in
standard position in $V_1$.  The isotopy deforms $\Gamma$ to the graph
$\Gamma' = e'_1 \cup e'_2$, with a single vertex $v'$ on $T$.  Since
$e'_2$ is a longitude on $\bdd V_1$ and $e'_1$ is in standard
position, $e'_1\cup e'_2$ bounds a $v'$-disk $\Delta$ in $V_1$.
Pushing $\Delta - v'$ to the interior of $V_1$ deforms $\Gamma'$ to a
graph in standard position, hence the result follows.

Now assume that $j(\a,\r) = q > 1$.  Choose a meridional disk $D_2$ of
$V_2$ containing $\a'$, intersecting $\r$ at $p$ points.  Since $\r'$
is a $(p,q)$ curve, and the jumping number of $\a$ is $q$, we can
choose the arc $\b$ on $\r'$ with $\bdd \b = \bdd \a$ so that the
interior of $\b$ is disjoint from $\bdd D_2$, hence $e''_2 = \a' \cup
\b$ is a longitude of $V_2$.  By an isotopy of $\Gamma\cap V_1$ we can
deform $e_2\cap V_1$ to $\b$, and $e_1$ to a loop $e'_1$ in standard
position in $V_1$.  Let $v' = e'_1 \cap e''_2$.  By an isotopy rel
$v'$ we can deform $e''_2$ to an edge $e'_2$ in $V_2$, which by
definition is in standard position in $V_2$ because $e''_2$ is a
longitude of $V_2$.  It follows that $\Gamma$ is isotopic to the graph
$\Gamma' = e'_1 \cup e'_2$ in standard position, hence $\Gamma$ is
standard.  \endproof

Suppose $F$ is a surface on the boundary of a 3-manifold $X$, and $c$
a simple closed curve in $F$.  Denote by $X_{c}$ the manifold obtained
from $X$ by attaching a 2-handle to $X$ along $c$, and by $F_{c}$ the
corresponding surface in $X_{c}$.  More explicitly, $X_{c} = X \cup
_{\varphi} (D^2 \times I)$, where $\varphi$ identifies $\bdd D^2
\times I$ to a regular neighborhood $A$ of $c$ in $F$, and $F_{c} =
(F-A)\cup (D^2 \times \bdd I)$.  We need the following version of 
handle addition lemma.

\proclaim{Lemma 8} Let F be a surface on the boundary of a 3-manifold
$X$, and $K$ a 1-manifold in $F$ with $F-K$ compressible in $X$.  Let
$c$ be a simple loop in $F-K$.  If $F_{c}$ has a compressing disk
$\Delta$ in $X_{c}$, then $F-c$ has a compressing disk $\Delta'$ in
$X$ such that $\bdd \Delta' \cap K \subset \bdd \Delta \cap K$.
\endproclaim

\proof This was proved in [Wu1].  Theorem 1 of [Wu1] says that under
the assumption of the lemma we have $|\bdd \Delta' \cap K| \leq |\bdd
\Delta \cap K |$, but that was proved by showing that $\bdd \Delta'
\cap K \subset \bdd \Delta \cap K$.  Note that when $K = \emptyset$,
it reduces to Jaco's Handle Addition Lemma [J, Lemma 1].  \endproof

\noindent
\demo{Proof of Theorem 3} By Lemma 4 we need only show that if the
exterior of any nontrivial subgraph of $\Gamma$ is a handlebody then
$\Gamma$ is standard.  Put $W = M - \eta(v)$, $t = W \cap \Gamma$, $X
= M -\eta(\Gamma) = W - \eta(t)$, and $P = \bdd W \cap X$.  By Lemma 5
we may assume that $P$ is incompressible, so by Lemma 2 and [MZ, Lemma
1], the manifold $M$ is a lens space $L(p,q)$.  Up to homeomorphism we
may assume $1\leq q \leq p/2$.

By Lemma 7 we may assume that $n\geq 3$, and by induction we may
assume that any nontrivial proper subgraph of $\Gamma$ is standard.
In particular, each $e_i$ is standard in $M$, so it is isotopic to a
core of either $V_1$ or $V_2$.  Since $n\geq 3$, at least two of the
$e_i$ are cores of the same $V_j$, hence up to relabeling we may
assume without loss of generality that $e_1$ and $e_2$ are both
isotopic to a core of $V_2$.

Consider the graph $\Gamma' = e_1 \cup ... \cup e_{n-1}$.  By
induction $\Gamma'$ is standard, so the edges are contained in two
$v$-disks if $M = L(p,q)$ with $2\leq q < p/2$, and one $v$-disk
otherwise.  Notice that in the first case the core of $V_1$ is
homotopic to $q$ times the core of $V_2$, so they represents different
elements in $\pi_1 M$.  Since by assumption $e_1$ and $e_2$ are
isotopic to the core of $V_2$, it follows that they are on the same
$v$-disk.  In either case there is a $v$-disk $D_1$ containing both
$e_1$ and $e_2$.  Taking a subdisk bounded by $e_1 \cup e_2$ and
pushing its interior off $D_1$, we get a $v$-disk $D_2$ bounded by
$e_1 \cup e_2$ with interior disjoint from $\Gamma'$.  Note that $D_2$
may intersect $e_n$.  However, the following sublemma says that $D_2$
can be rechosen to have interior disjoint from $e_n$ as well.

\proclaim{Sublemma} There is a $v$-disk $D_3$ bounded by $e_1\cup e_2$
with interior disjoint from $\Gamma$.  \endproclaim

\proof Consider the handlebody $X = M- \eta(\Gamma)$.  Let $c_i$ be
the meridian curve of $e_i$ on $F = \bdd X$, and put ${C} =
\{c_1, ... , c_n\}$.  Let $K = c_1 \cup ... \cup c_{n-1}$.  By Lemma
2, the tangle $(W, t)$ is completely tubing compressible, so $K$ is a
primitive set on $\bdd X$, hence $F-K$ is compressible.  We now apply
Lemma 8 to $(X, F, K, c)$ with $c = c_n$.  Note that after attaching a
2-handle to $c_n$, the manifold $X' = X_{c_n}$ is the same as the
exterior of the graph $\Gamma' = e_1 \cup ... \cup e_{n-1}$, and the
surface $F_{c_n} = \bdd X'$.

Recall that $e_1 \cup e_2$ bounds a $v$-disk $D_2$ in $M$ with
interior disjoint from $\Gamma'$, so its restriction to $X' = X_{c_n}$
is a compressing disk $\Delta$ of $\bdd X' = F_{c_n}$ intersecting
each of $c_1$ and $c_2$ at a single point, and is disjoint from $c_3,
..., c_{n-1}$.  Therefore, by Lemma 8, there is a compressing disk
$\Delta'$ of $F-c_n$ in $X$, such that $\bdd \Delta'$ intersects each
of $c_1$ and $c_2$ at most once, and is disjoint from $c_3, ...,
c_{n-1}$.  Since it is a compressing disk of $F - c_n$, it is also
disjoint from $c_n$.

Now $\bdd \Delta'$ cannot be disjoint from $C$, because we have
assumed that the surface $P$ homotopic to $F - C$ is
incompressible.  Also, $\bdd \Delta' \cap C$ cannot be a single
point in $c_1$, say, because then the frontier of a regular
neighborhood of $\Delta' \cup c_1$ would be a compressing disk of $F -
C$, which is again a contradiction.  It follows that $\bdd
\Delta'$ intersects each of $c_1$ and $c_2$ at exactly one point, and
is disjoint from the other $c_j$'s.  Since $\Gamma$ is a spine of
$N(\Gamma)$, by shrinking $N(\Gamma)$ to $\Gamma$, the disk $\Delta'$
becomes a $v$-disk $D_3$ in $M$ bounded by $e_1 \cup e_2$, with
interior disjoint from $\Gamma$.  This completes the proof of the
sublemma.  \endproof

We now continue to show that $\Gamma$ is standard in $M$.  By
induction we may assume that $\Gamma'' = e_2 \cup ... \cup e_n$ is in
standard position in $M = V_1 \cup V_2$, with $e_2$ on a $v$-disk $D'$
in $V_2$, say, which contains all the edges of $\Gamma''$ in $V_2$.
Consider the disk $D_3$ bounded by $e_1\cup e_2$ as given by the
sublemma.  It has interior disjoint from $\Gamma$, so by considering
$D_3 \cap D'$ and using an innermost circle outermost arc argument one
can show that $D_3$ can be modified so that it intersects $D'$ only
along the edge $e_2$.  Pushing the part of $D_3$ near $e_2$ slightly
off $e_2$, we get a $v$-disk $D_4$ with boundary the union of $e_1$
and a loop $e'_1$ on $D'$, which is a parallel copy of $e_2$
intersecting $\Gamma$ only at $v$.  One can then isotope $e_1$ via the
disk $D_4$ to the edge $e_1'$, which lies on the $v$-disk $D'$.  Thus
after this isotopy all edges of $\Gamma$ are now contained in the
$v$-disks which contain $\Gamma''$.  Therefore $\Gamma$ is also
standard by definition.  \endproof

\noindent
\demo{Proof of Theorem 1} Suppose $(W,t)$ is completely tubing
compressible.  Then by Lemma 2 the corresponding graph $\hat t$ in
$\hat W = W \cup B$ has the property that the exterior of any proper
subgraph of $\hat t$ is a handlebody.  By Theorem 3, $\hat t$ is
contained in the union of at most two $v$-disks $D_1$ and $D_2$, with
$D_i$ in $V_i$.  By an isotopy rel $\hat t$ we may assume that $D_i
\cap \bdd W$ consists of two arcs, hence $D_1 \cap W$ and $D_2 \cap W$
are two disjoint bands in $W$ containing $t$, and the result follows.
\endproof

\enddocument